\documentclass{ws-ijaitarxiv}
\makeatletter
\def\blfootnote{\gdef\@thefnmark{}\@footnotetext}
\makeatother

\usepackage{tikz}
\usetikzlibrary{arrows,positioning,shapes.geometric}

\usepackage{pgfplots}
\pgfplotsset{compat=1.8}

\newcommand{\R}{\mathbb{R}}

\DeclareMathOperator{\Ima}{Im}
\DeclareMathOperator{\rank}{rank}
\DeclareMathOperator{\Dom}{Dom}

\begin{document}

\markboth{C.\ Wu \& C.\ A.\ Hargreaves}{Topological Machine Learning for Mixed Numeric and Categorical Data}

%
\catchline{}{}{}{}{}
%

\title{Topological Machine Learning for Mixed Numeric and Categorical Data}

\author{Chengyuan Wu\textsuperscript{*,$\dagger$} \and Carol Anne Hargreaves\textsuperscript{*,$\ddagger$}}

\address{\textsuperscript{*}Data Analytics Consulting Centre, Department of Statistics and Applied Probability\\ 
National University of Singapore, Singapore 117546, Singapore\\
\email{\textsuperscript{$\dagger$}stawuc@nus.edu.sg}
\email{\textsuperscript{$\ddagger$}stacah@nus.edu.sg}
\blfootnote{\textsuperscript{$\dagger$}Corresponding author.}
}

\maketitle

\begin{history}
\received{(Day Month Year)}
\revised{(Day Month Year)}
\accepted{(Day Month Year)}
\end{history}

\maketitle


\begin{abstract}
Topological data analysis is a relatively new branch of machine learning that excels in studying high-dimensional data, and is theoretically known to be robust against noise. Meanwhile, data objects with mixed numeric and categorical attributes are ubiquitous in real-world applications. However, topological methods are usually applied to point cloud data, and to the best of our knowledge there is no available framework for the classification of mixed data using topological methods. In this paper, we propose a novel topological machine learning method for mixed data classification. In the proposed method, we use theory from topological data analysis such as persistent homology, persistence diagrams and Wasserstein distance to study mixed data. The performance of the proposed method is demonstrated by experiments on a real-world heart disease dataset. Experimental results show that our topological method outperforms several state-of-the-art algorithms in the prediction of heart disease.
\end{abstract}


\section{Introduction}
Topological data analysis (TDA) is a relatively new subject that is gaining popularity in many fields, such as network analysis,\cite{Carstens2013,horak2009persistent} biomolecular chemistry,\cite{xia2018multiscale,xia2015multiresolution} and drug design.\cite{cang2018integration,wu2018topp} Topological data analysis is often referred to as studying the ``shape'' of data, in order to deduce fundamental characteristics of the data. The primary tool used in TDA is persistent homology,\cite{edelsbrunner2012persistent,Zomorodian2005} though there are also other tools such as Mapper,\cite{nicolau2011topology,singh2007topological} discrete Morse theory, \cite{forman2002user,reininghaus2010tadd,wu2020discrete} as well as other techniques from algebraic topology.\cite{hansen2019toward,letscher2012persistent,wu2020weighted,wu2019magnus} It is generally acknowledged that topological data analysis is effective at analyzing high-dimensional noisy data.\cite{adams2017persistence,offroy2016topological} We also remark that topological methods have also recently gained prominence in physics, with the 2016 Nobel Prize in Physics being awarded for theoretical discoveries of topological phase transitions and topological phases of matter.\cite{haldane2017nobel}

In real-world applications, data sets often have both numeric and categorical attributes. The coexistence of numeric and categorical variables often makes machine learning methods designed for single-type data inapplicable to mixed-type data.\cite{hsu2008extended,ji2015initialization} Traditionally, TDA is usually applied to point cloud data or spatial data.\cite{Ghrist2008} The strengths of TDA include the property of being coordinate-free\cite{lum2013extracting,offroy2016topological} (independent of the coordinate system chosen), as well as being translation-invariant and rotation-invariant. \cite{bonis2016persistence,khasawneh2016chatter} A drawback of these strengths is that it may be hard for TDA to effectively analyze data that is sensitive to choice of coordinates, translation, and/or rotation. Examples of such data include data with heterogeneous features, where each coordinate represents a fundamentally different feature (e.g.\ light, temperature, humidity).\cite{wu2019topological} To the best of our knowledge, there is currently no readily available framework for the classification of mixed numeric and categorical data using TDA. In view of the ubiquity of mixed-type data and the rising popularity of TDA, it is of interest to develop a topological machine learning method for mixed data.

\begin{figure}[htbp]
\begin{center}
\begin{tikzpicture}[>=latex']
        \tikzset{block/.style= {draw, rectangle, align=center,minimum width=2cm,minimum height=1cm},
        rblock/.style={draw, shape=rectangle,rounded corners=1.5em,align=center,minimum width=2cm,minimum height=1cm},
        input/.style={shape=ellipse, aspect=2.2, draw},
        }
        \node [rblock]  (a1) {Mixed data};
        \node [block, right =2em of a1] (a2) {One-hot \\ encoding};
        \node [block, right =2em of a2] (a3) {Standardization};
        \node [block, right =2em of a3] (a4) {Symmetry \\ breaking};
        \node [block, below right =3em and -5em of a1] (b1) {Projection \\ maps};
        \node [block, right =1em of b1] (b2) {Point \\ clouds};
        \node [block, right =1em of b2] (b3) {Persistence \\ diagrams};
        \node [block, right =1em of b3] (b4) {$k$-NN \\ (Wasserstein \\ distance)};
        \node [input, right =1em of b4] (b5) {Classification};
        \node [coordinate, below right =1em and 3em of a4] (right) {};  
        \node [coordinate,above left =2em and 1em of b1] (left) {};  

        \path[draw,->] (a1) edge (a2)
                    (a2) edge (a3)
                    (a3) edge (a4)
                    (a4.east) -| (right) -- (left) |- (b1)
                    (b1) edge (b2)
                    (b2) edge (b3)
                    (b3) edge (b4)
                    (b4) edge (b5)
                    ;
    \end{tikzpicture}
    \end{center}
\caption{Basic workflow of Topological Machine Learning for Mixed Numeric and Categorical Data.}
\label{fig:workflow}
\end{figure}
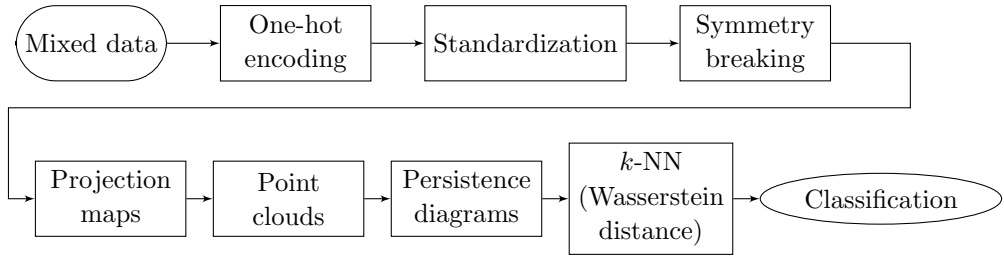

In this paper, we propose a novel topological machine learning method for mixed data (TopMix). In our method, the categorical variables are first converted to binary variables via one-hot encoding. All predictor variables are subsequently standardized, and a basic \emph{symmetry breaking}\cite{wu2019topological} technique is applied to the data for TDA to better deal with heterogeneous features. Subsequently, each data point is converted into a point cloud via multiple projection maps. We then generate persistence diagrams from the point cloud data, and calculate the Wasserstein distance between the persistence diagrams. Lastly, we use the $k$-nearest neighbors algorithm ($k$-NN) for supervised machine learning (classification). The basic workflow of our paper is summarized in Figure \ref{fig:workflow}.

We remark that the technique of converting each data point into a point cloud via multiple projection maps is specialized for the setting of mixed numeric and categorical data. Hence, it is a new innovation that is not present in the authors' previous paper\cite{wu2019topological}, which focuses on the setting of multivariate time series data.

For applications, we apply our method to heart disease prediction. Heart disease is the leading cause of death in the industrialized world.\cite{lauer1999cause} For instance, in 2002, 696,947 people in the United States died of heart disease, compared with 557,271 deaths from cancer.\cite{twombly2005cancer} We use a dataset originating from the seminal paper by R.\ Detrano et al.\cite{detrano1989international} In the dataset, there are 14 attributes including numeric and categorical variables. The goal is to predict whether a patient has heart disease ($>50\%$ luminal narrowing of any major epicardial vessel) or not. We show that topological methods are effective in predicting heart disease using mixed data. Our topological method outperforms several state-of-the-art algorithms in the classification of heart disease.

The rest of the paper is organized as follows. We first review some related work in Section \ref{sec:relatedwork}. This is followed by a brief introduction to the background information on TDA in Section \ref{sec:background}. In Section \ref{sec:tml}, we present our topological machine learning method for mixed numeric and categorical data. In Section \ref{sec:experiment}, we report the experimental results, which demonstrate the viability of the proposed method. Finally, we draw conclusions in Section \ref{sec:conclusion}. 

\subsection{Related Work}
\label{sec:relatedwork}
In the paper by X.\ Ni et al.,\cite{ni2017composing} the authors proposed a clustering method for mixed data based on a tree-structured graphical model. Their tree-structured model factorizes into a product of pairwise interactions. Furthermore, the authors leverage theory from TDA to adaptively merge trivial peaks of the density function into larger ones in order to achieve meaningful clusterings. Persistent homology theory is used to automatically determine the number of clusters in the data. An earlier seminal paper by Chazal et al.\cite{chazal2013persistence} introduced the novel idea of using topological persistence to guide the merging of clusters. Their algorithm provides additional feedback in the form of a persistence diagram, which the authors prove to reflect the prominences of the modes of the density. The algorithm requires rough estimates of the density at the data points, and knowledge of approximate pairwise distances between them, and hence is applicable in any metric space. Their method can be theoretically proven to output the correct number of clusters under certain mild sampling conditions.

In recent years, topological techniques have been effectively combined with machine learning or statistical methods. In the paper by C.\ Hofer et al.,\cite{hofer2017deep} the authors introduced a technique that enables the input of topological signatures to deep neural networks for learning a task-optimal representation during training. An advantage of their method is that it learns the representation instead of mapping topological signatures to a pre-defined representation. P.\ Bubenik defined the persistence landscape,\cite{bubenik2015statistical}  which is a novel topological summary for data. Since this summary lies in a vector space, it is possible to combine it with tools from statistics and machine learning.  A number of standard statistical tests can be used for statistical inference using persistence landscapes, for example the two-sample $Z$-test and Hotelling's $T^2$ test.

C.\ Wu and C.\ A.\ Hargreaves\cite{wu2019topological} developed a framework for analyzing multivariate time series using TDA. The methodology includes converting the multivariate time series to point cloud data, calculating Wasserstein distances between the persistence diagrams, and using the $k$-NN algorithm for classification. For applications, the authors focus on room occupancy detection based on 5 time-dependent variables (temperature, humidity, light, CO\textsubscript{2} and humidity ratio).

In the paper by J.\ Ji et al.,\cite{ji2015initialization} the authors proposed a new initialization method for mixed data clustering. Prior to their paper, most of the initialization approaches are dedicated to partitional clustering algorithms which process either categorical or numerical data only. In the paper, the authors introduced a new definition of density to assess the cohesiveness of data objects with mixed numeric and categorical attributes.

A.\ Ahmad and L.\ Dey\cite{ahmad2007k} presented a clustering algorithm that works well for data with mixed numeric and categorical features. The authors proposed a new cost function and distance measure based on co-occurrence of values. In their scheme, $\delta (p,q)$ which denotes the distance between a pair of distinct values $p$ and $q$ of an attribute, is computed as a function of their co-occurrence with other attribute values. The contribution of a categorical attribute is inherent in the distance measure itself and need not be user defined.

The paper by J.\ Nahar et al.\cite{nahar2013computational} investigates various computational intelligence techniques in the detection of heart disease. In the paper, the Cleveland dataset\cite{detrano1989international} from the UCI Machine Learning Repository\cite{Dua:2019} is used. In particular, the authors highlight the potential of a medical knowledge driven feature selection process for heart disease diagnosis. Experiments show that the medical knowledge based feature selection method has shown promise for use in heart disease diagnostics.

R.\ Das, I.\ Turkoglu and A.\ Sengur\cite{das2009effective} explored the effective diagnosis of heart disease through neural network ensembles. Ensemble based methods can enable an increase in performance by combining several individual neural networks to train on the same task. The authors utilize SAS base software 9.1.3 in their methodology, and achieved good results using three independent neural network models in the ensemble model.

\section{Background}
\label{sec:background}

We give a brief overview of the key concepts in TDA and persistent homology, and refer the reader to the appropriate references for more details. A classical text for algebraic topology is the book by A.\ Hatcher.\cite{Hatcher2002} The survey article by H.\ Edelsbrunner and J.\ Harer,\cite{edelsbrunner2008persistent} as well as the review paper by R.\ Ghrist,\cite{Ghrist2008} provide a superb introduction to persistent homology. In addition, the paper by A.\ Zomorodian and G.\ Carlsson\cite{Zomorodian2005} gives a comprehensive overview of persistent homology from a mathematical and computational perspective.

\subsection{Simplicial complexes}
Simplicial complexes are one of the main objects of study in algebraic topology. A simplicial complex can be regarded as a set composed of vertices, edges, triangles, and higher dimensional simplices.

\begin{figure}[!htbp]
\begin{center}
\begin{tikzpicture}
\node[circle,fill=black,inner sep=0pt,minimum size=3pt] (a) at (0,0){};
\node[below] at (0,-0.1) {$v_0$};
\node[below] at (0,-1) {vertex $\{v_0\}$};
\end{tikzpicture}
\begin{tikzpicture}
\node[circle,fill=black,inner sep=0pt,minimum size=3pt] (a) at (0,0){};
\node[circle,fill=black,inner sep=0pt,minimum size=3pt] (b) at (1.5,0){};
\node[below] at (0,-0.1) {$v_0$};
\node[below] at (1.5,-0.1) {$v_1$};
\node[below] at (0.75,-1) {edge $\{v_0,v_1\}$};
\draw (a) -- (b);
\end{tikzpicture}
\begin{tikzpicture}
\coordinate (a) at (0,0);
\coordinate (b) at (1.5,0);
\coordinate (c) at (0.75,1.3);
\filldraw[gray!50] (a.center) -- (b.center) -- (c.center) -- cycle;
\draw (a.center) -- (b.center) -- (c.center) -- cycle;
\node[circle,fill=black,inner sep=0pt,minimum size=3pt] at (a){};
\node[circle,fill=black,inner sep=0pt,minimum size=3pt] at (b){};
\node[circle,fill=black,inner sep=0pt,minimum size=3pt] at (c){};
\node[below] at (0,-0.1) {$v_0$};
\node[below] at (1.5,-0.1) {$v_1$};
\node[left] at (0.65,1.3) {$v_2$};
\node[below] at (0.75,-1) {triangle $\{v_0,v_1,v_2\}$};
\end{tikzpicture}
\begin{tikzpicture}
\coordinate (a) at (0,0);
\coordinate (b) at (1.5,0);
\coordinate (c) at (0.75,1.3);
\coordinate (d) at (1.8,0.6);
\begin{scope}[dashed,,opacity=0.5]
\draw (a) -- (d);
\end{scope}
\draw[fill=gray,opacity=0.5] (a)--(b)--(c);
\draw[fill=gray,opacity=0.5] (b)--(c)--(d);
\draw (a.center) -- (b.center) -- (c.center) -- cycle;
\draw (c) -- (d);
\draw (b) -- (d);
\node[circle,fill=black,inner sep=0pt,minimum size=3pt] at (a){};
\node[circle,fill=black,inner sep=0pt,minimum size=3pt] at (b){};
\node[circle,fill=black,inner sep=0pt,minimum size=3pt] at (c){};
\node[circle,fill=black,inner sep=0pt,minimum size=3pt] at (d){};
\node[below] at (0,-0.1) {$v_0$};
\node[below] at (1.5,-0.1) {$v_1$};
\node[left] at (0.65,1.3) {$v_2$};
\node[right] at (1.8,0.7) {$v_3$};
\node[below] at (0.5,-1) {tetrahedron $\{v_0,v_1,v_2,v_3\}$};
\end{tikzpicture}
\caption{A 0-simplex (vertex), 1-simplex (edge), 2-simplex (triangle) and 3-simplex (tetrahedron).}
\label{fig:simplices}
\end{center}
\end{figure}
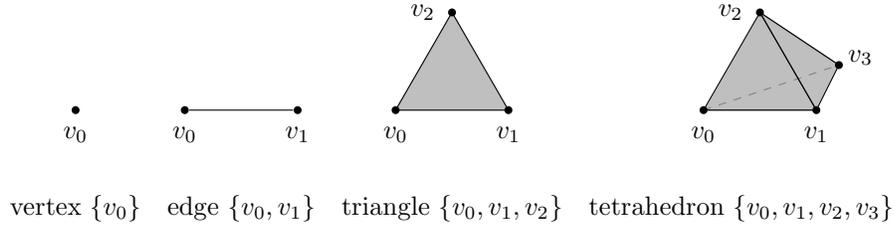

More formally, a \emph{simplicial complex} $K$ is a collection of sets such that $\sigma\in K$ and $\tau\subseteq\sigma$ implies $\tau\in K$. The sets $\sigma\in K$ are called the \emph{simplices} of the simplicial complex $K$. We call the singleton sets $\{v\}$ the \emph{vertices} of $K$. The dimension of a simplex $\sigma\in K$ is defined to be $\dim(\sigma)=|\sigma|-1$, and we call a simplex of dimension $k$ a \emph{$k$-simplex}. Simplices of dimension 0, 1, 2, 3 represent a \emph{vertex}, \emph{edge}, \emph{triangle} and \emph{tetrahedron} respectively, as shown in Figure \ref{fig:simplices}.

A type of simplicial complex frequently used in TDA is the \emph{Vietoris-Rips complex} (or \emph{Rips complex} for short).

\begin{definition}
Let $\{x_i\}$ be a set of points in the Euclidean space $\mathbb{R}^n$. The Rips complex $\mathcal{R}_\epsilon$ is the simplicial complex whose $k$-simplices consist of each subset of $k+1$ points $\{x_j\}_{j=0}^k$ which are pairwise within distance $\epsilon$.
\end{definition}

\begin{definition}
Let $K$ be a simplicial complex. Suppose $L$ is a simplicial complex such that every face of $L$ belongs to $K$, that is, $L\subseteq K$. We say that $L$ is a \emph{simplicial subcomplex} of $K$.
\end{definition}

We also introduce the notion of a \emph{filtration} of a simplicial complex $K$, which is a nested sequence of complexes $\emptyset=K^0\subseteq K^1\subseteq\dots\subseteq K^m=K$. We say that $K$ is a \emph{filtered complex}.

\subsection{Homology}
The $k$th \emph{chain group} $C_k$ of a simplicial complex $K$ is defined to be the free abelian group with basis to be the set of oriented $k$-simplices. The boundary operator $\partial_k: C_k\to C_{k-1}$ is defined on an oriented simplex $\sigma=[v_0,v_1,\dots,v_k]$ by
\begin{equation*}
\partial_k(\sigma)=\sum_{i=0}^k (-1)^i [v_0,\dots,\hat{v_i},\dots,v_k],
\end{equation*}
where $\hat{v_i}$ denotes the deletion of the vertex $v_i$.

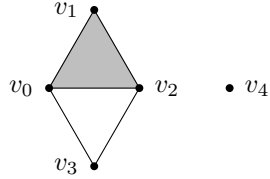
\begin{figure}[htbp]
\begin{center}
\begin{tikzpicture}[scale=0.8]
\coordinate (a) at (0,0);
\coordinate (b) at (1.5,0);
\coordinate (c) at (0.75,1.3);
\coordinate (d) at (0.75,-1.3);
\coordinate (e) at (3,0);
\filldraw[gray!50] (a.center) -- (b.center) -- (c.center) -- cycle;
\draw (a.center) -- (b.center) -- (c.center) -- cycle;
\draw (a) --(d) --(b);
\node[circle,fill=black,inner sep=0pt,minimum size=3pt] at (a){};
\node[circle,fill=black,inner sep=0pt,minimum size=3pt] at (b){};
\node[circle,fill=black,inner sep=0pt,minimum size=3pt] at (c){};
\node[circle,fill=black,inner sep=0pt,minimum size=3pt] at (d){};
\node[circle,fill=black,inner sep=0pt,minimum size=3pt] at (e){};
\node[left] at (-0.1,0) {$v_0$};
\node[right] at (1.6,0) {$v_2$};
\node[left] at (0.65,1.3) {$v_1$};
\node[left] at (0.65,-1.3) {$v_3$};
\node[right] at (3.1,0) {$v_4$};
\end{tikzpicture}
\caption{For the above simplicial complex, we have $\beta_0=2$ (2 connected components), $\beta_1=1$ (1 ``circular'' hole which corresponds to the unshaded region) and $\beta_2=0$ (no ``voids'').}
\label{fig:betti}
\end{center}
\end{figure}

Subsequently, the $k$th \emph{homology group} is defined as the quotient $H_k=Z_k/B_k$, where $Z_k=\ker\partial_k$ and $B_k=\Ima\partial_{k+1}$ are the \emph{cycle group} and the \emph{boundary group} respectively. The rank of the $k$th homology group $\beta_k=\rank(H_k)$ can be said to count the number of $k$-dimensional ``holes'' in $K$, as illustrated in Figure \ref{fig:betti}.

\subsection{Persistent homology}
Given a filtered complex $K$, we may define the corresponding boundary operators $\partial_k^i$ and groups $C_k^i$, $Z_k^i$, $B_k^i$ and $H_k^i$ for the $i$th complex $K^i$. The \emph{$p$-persistent $k$th homology group} of $K^i$ is defined as
\begin{equation*}
H_k^{i,p}=Z_k^i/(B_k^{i+p}\cap Z_k^i).
\end{equation*}

The filtered complex $K$ is usually obtained by the construction of Rips complexes over a range of distances $\epsilon$. Persistent homology detects those topological features which persist over a parameter range, revealing meaningful structures in the data.

\section{Topological Machine Learning for Mixed Numeric and Categorical Data}
\label{sec:tml}

In this section, we describe our approach of using topological machine learning methods to analyze mixed data. A basic summary of the workflow can be found in Figure \ref{fig:workflow}.

\subsection{Notation}
We first introduce a standard notation for mixed data, following the paper by Z.\ Huang and M.\ K.\ Ng,\cite{huang1999fuzzy} as well as J.\ Ji et al.\cite{ji2015initialization} Let $\mathbf{X}=\{X_1,X_2,\dots,X_n\}$ denote a dataset of $n$ data objects. Each object $X_i$ has $m$ attributes $A_1,A_2,\dots,A_m$. We represent each $X_i$ as a $m$-tuple $(x_{i,1},x_{i,2},\dots,x_{i,m})$. Each attribute $A_j$ is associated with a domain of values, denoted by $\Dom(A_j)$, which is either numeric (real numbers) or categorical (finite, unordered set). A categorical domain is generally represented by $\Dom(A_j)=\{a_{j,1},a_{j,2},\dots,a_{j,s}\}$, where $s$ is the number of possible categorical values for the categorical attribute $A_j$. Each data object $X_i$ can be logically represented as a conjunction of attribute-value pairs: 
\begin{equation*}
[A_1=x_{i,1}]\wedge [A_2=x_{i,2}]\wedge\dots\wedge [A_m=x_{i,m}].
\end{equation*}

\subsection{Methodology}
\subsubsection{One-hot encoding}
Firstly, we apply one-hot encoding to the mixed data, converting each categorical variable with $s$ possible values to $s$ binary variables. That is, each categorical attribute $A_j$ with $\Dom(A_j)=\{a_{j,1},a_{j,2},\dots,a_{j,s}\}$ is replaced with $s$ binary attributes $B_1,B_2,\dots,B_s$, with $\Dom(B_i)=\{0,1\}$ for $1\leq i\leq s$.

\subsubsection{Standardization}
We standardize all variables (including binary variables) to have zero mean and unit variance. This is to ensure that all variables are on the same scale, preventing a feature with larger scale from dominating other features. We remark that the standardization of binary variables is also done in the algorithms KNNImpute (for categorical data)\cite{stekhoven2012missforest} and Lasso.\cite{tibshirani1997lasso}

\subsubsection{Symmetry breaking}
\label{sec:symmbreak}
Symmetry breaking refers to adding a fixed constant vector to each data object, with the purpose of enabling TDA methods to better distinguish point clouds that may just differ by translation or rotation. Symmetry breaking was introduced in the context of studying multivariate time series using topological methods.\cite{wu2019topological} Basically, symmetry breaking attempts to ``disable'' the translational / rotational invariance property of TDA for data that do not require it. 

\begin{definition}
Let $X=(x_1,x_2,\dots,x_m)$ be a data object represented as a $m$-tuple in $\R^m$. Let $\mathbf{v}=(c_1,c_2,\dots,c_m)$ be a fixed vector in $\R^m$. We define the new data object $X'$ obtained by \emph{symmetry breaking} (of $X$) to be $X'=X+\mathbf{v}$.
\end{definition}

An example of the fixed vector is $\mathbf{v}=(5,6,7,\dots,m+4)$. We will be using this fixed vector in the paper. In Section \ref{sec:explainsymm}, we will illustrate how symmetry breaking can be useful in analyzing data with heterogeneous features, as well as explain our heuristic choice of fixed vector $\mathbf{v}$.

\subsubsection{Projection maps}
After applying symmetry breaking, the new data object
\begin{equation*}
X'=(x_1+c_1,x_2+c_2,\dots,x_m+c_m)
\end{equation*}
is a single point in $\R^m$. However, a single point has trivial homology and trivial persistent homology, hence we will need a point cloud (set of multiple data points in Euclidean space) in order for topological methods to work. In contrast, in the authors' previous paper on multivariate time series\cite{wu2019topological}, the point cloud resulting from a time window of length $w>1$ already consists of multiple data points. Hence, there was no need for the technique of projection maps in the setting of multivariate time series data.

\begin{definition}
The $i$th \emph{projection map} $p_i:\R^m\to\R^m$ is defined by
\begin{equation*}
p_i(y_1,y_2,\dots,y_i,\dots,y_m)=(y_1,y_2,\dots,0,\dots,y_m).
\end{equation*}
\end{definition}

The projection map $p_i$ changes the $i$th coordinate of a vector to 0. This is equivalent to projecting the vector onto the hyperplane $H=\{(x_1,x_2,\dots,x_m)\in\R^m\mid x_i=0\}$. The projection map defined above is idempotent, namely $p_i\circ p_i=p_i$.

\subsubsection{Point clouds}
\label{sec:pointcloud}
We define the point cloud $S(X')$ associated to the data object $X'\in\R^m$ to be
\begin{equation*}
S(X')=\{X'\}\cup\{p_1(X'),p_2(X'),\dots,p_m(X')\}.
\end{equation*}

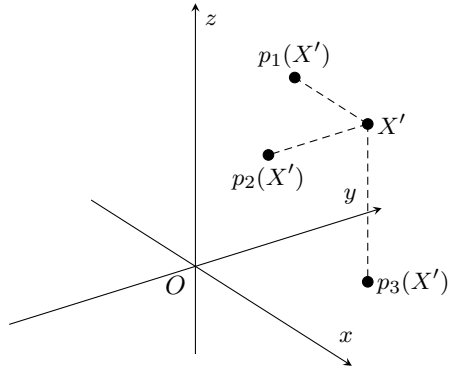
\begin{figure}[htbp]
\begin{center}
\begin{tikzpicture}
\begin{axis}[
  ticks=none,
  view={55}{30},
  axis lines=center,
  width=10cm,height=10cm,
  xmin=-10,xmax=15,ymin=-15,ymax=15,zmin=-5,zmax=15,
  xlabel={$x$},ylabel={$y$},zlabel={$z$},
  label style={font=\small}
]

\addplot3 [only marks] coordinates {(7,8,9) (0,8,9) (7,0,9) (7,8,0)};

\addplot3 [no marks,densely dashed] coordinates {(7,8,9) (0,8,9)};
\addplot3 [no marks,densely dashed] coordinates {(7,8,9) (7,0,9)};
\addplot3 [no marks,densely dashed] coordinates {(7,8,9) (7,8,0)};

\node [right] at (axis cs:7,8,9) {\small $X'$};
\node [above] at (axis cs:0,8,9) {\small $p_1(X')$};
\node [below] at (axis cs:7,0,9) {\small $p_2(X')$};
\node [right] at (axis cs:7,8,0) {\small $p_3(X')$};
\node [below left] at (axis cs:0,0,0) {$O$};
\end{axis}
\end{tikzpicture}
\caption{For $m=3$, the point cloud $S(X')$ consists of $X'$, as well as the 3 projected points $p_1(X')$, $p_2(X')$ and $p_3(X'$). The 3 projected points are projections of $X'$ onto the $yz$-, $xz$-, and $xy$-planes respectively.}
\label{fig:projpointcloud}
\end{center}
\end{figure}

That is, $S(X')$ consists of $m+1$ points, namely the point $X'$ as well as the $m$ projected points $p_1(X'),p_2(X'),\dots,p_m(X')$. We show an example for the case $m=3$ in Figure \ref{fig:projpointcloud}.

The point cloud $S(X')$ contains intrinsic information about the data object $X'$ in the form of distances between the points in $S(X')$. For instance, if $X'=(y_1,y_2,\dots,y_m)$, we can calculate the following Euclidean distances between $X'$ and its projections:
\begin{align}
d(X',p_i(X'))&=|y_i|, \label{eq:dist1}\\
d(p_i(X'),p_j(X'))&=\sqrt{y_i^2+y_j^2}, \quad\text{for $i\neq j$}. \label{eq:dist2}
\end{align}

\subsubsection{Persistence diagrams}
A persistence diagram\cite{cohen2007stability} is a multiset of points in the space $\Delta:=\{(b,d)\in\R^2\mid b,d\geq 0, b\leq d\}$. Each point $(b,d)$ represents a persistent generator (of a given dimension), where $b$ denotes the birth of the generator and $d$ its death. In brief, the persistence diagram is a visual representation of the persistent homology of a point cloud. The persistence diagram is independent of choice of generators and hence is unique.\cite{cohen2010lipschitz} A notable result is the stability of persistence diagrams with respect to Hausdorff distance, bottleneck distance,\cite{cohen2007stability} as well as Wasserstein distance.\cite{cohen2010lipschitz} Such stability results give TDA the benefit of being robust to noise.

\begin{figure}[!htbp]
\begin{center}
\includegraphics[scale=0.7]{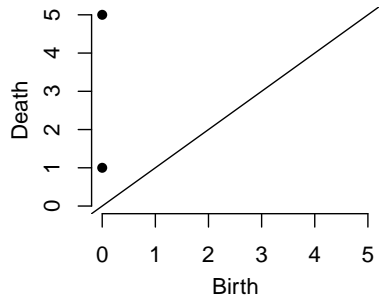}
\caption{The persistence diagram (dimension 0) for the point cloud $S$.}
\label{fig:pdexample}
\end{center}
\end{figure}

For the consideration of readability, we include a concrete example that illustrates the relationship between the persistent homology of a point cloud and its persistence diagram. Consider the point cloud $S=\{(0,0),(1,0)\}$ consisting of two points (0-simplices). At the start of the filtration process, there are two separate connected components (namely the two 0-simplices in $S$), hence this corresponds to two points on the persistence diagram (Figure \ref{fig:pdexample}) with birth time 0. At the filtration stage of $\epsilon=1$, the Rips complex $\mathcal{R}_\epsilon$ now consists of only one single connected component (namely the 1-simplex consisting of the two points in $S$ and the edge joining them). This corresponds to a death time of 1, which explains the point (0,1) on the persistence diagram. This 1-simplex theoretically persists to infinity (death time of infinity), but for practical purposes in the code we have to set a maximum value of $\epsilon$ for the Rips filtration (in this case $\texttt{maxscale}=5$). Hence, this explains the point (0,5) on the persistence diagram.

\subsubsection{$k$-NN (Wasserstein distance)}
The Wasserstein distance\cite{berwald2018computing,cohen2010lipschitz,mileyko2011probability} is commonly used to compare between two persistence diagrams.

\begin{definition}
The $p$-th Wasserstein distance between two persistence diagrams $D_1$, $D_2$ (of the same dimension) is defined to be
\begin{equation*}
W_p(D_1,D_2)=\left(\inf_{\varphi:D_1\to D_2}\sum_{x\in D_1}\| x-\varphi(x)\|^p_\infty\right)^{1/p},
\end{equation*}
where the infimum is taken over all bijections $\varphi$ between $D_1$ and $D_2$.
\end{definition}

As $p$ tends to infinity, the Wasserstein distance $W_p$ approaches the bottleneck distance $W_\infty$. The bottleneck distance captures the most perturbed topological feature (or the extreme behavior) of a point cloud, and can lead to noisier results than the Wasserstein distance.\cite{hajij2018visual} 

For this paper, we will use the Wasserstein distance with $p=1$, also known as the 1-Wasserstein distance or ``earth mover's distance''. The 1-Wasserstein distance is widely utilized in computer science,\cite{rabin2009statistical,rubner2000earth} including a recent usage in generative adversarial networks.\cite{arjovsky2017wasserstein}

Subsequently, to carry out classification (supervised machine learning), we use the $k$-nearest neighbors algorithm ($k$-NN) based on the Wasserstein distance. For each point cloud $S(X')$ (corresponding to a data object $X'$) in the test set, we will determine its $k$-nearest neighbors $\{S(Y_1),S(Y_2),\dots,S(Y_k)\}$ in the training set, with respect to the Wasserstein distance. Finally, we classify $X$ based on the majority class of the elements in the set $\{Y_1,Y_2,\dots,Y_k\}$.

\subsection{Elaboration on symmetry breaking}
\label{sec:explainsymm}
In this section, we illustrate how symmetry breaking, together with projection maps, can be helpful in analyzing data with heterogeneous features.

Consider two data objects $X=(1,2)$ and $Y(2,1)$. Their associated point clouds are $S(X)=\{(1,2),(0,2),(1,0)\}$ and $S(Y)=\{(2,1),(0,1),(2,0)\}$. We note that the pairwise distances between points in $S(X)$ are exactly the same as the respective pairwise distances between points in $S(Y)$, namely $1$, $2$ and $\sqrt{5}$. This would mean that topological methods will not be able to distinguish between $S(X)$ and $S(Y)$. The basic principle is that topological methods does not distinguish between point clouds that are related by ``symmetry'' (e.g.\ differ by rotation, translation, reflection).

Now, consider $\mathbf{v}=(5,6)$ such that we have
\[X'=X+\mathbf{v}=(6,8)\]
and
\[Y'=Y+\mathbf{v}=(7,7).\]
Then, the associated point clouds become $S(X')=\{(6,8),(0,8),(6,0)\}$ and $S(Y')=\{(7,7),(0,7),(7,0)\}$. The pairwise distances between points in $S(X')$ are $6$, $8$ and $10$, while the pairwise distances between points in $S(Y')$ are $7$, $7$ and $7\sqrt{2}$. Due to the difference in distances, TDA will be able to tell apart the point clouds $S(X')$ and $S(Y')$, which is the desired outcome.

Next, we will explain our heuristic choice of fixed vector $\mathbf{v}=(5,6,7,\dots,m+4)$ as mentioned in Section \ref{sec:symmbreak}. The main reason is to try to make the components $y_i$ in the data object $X'=(y_1,y_2,\dots,y_m)$ all positive (or mostly positive). By observing Equations \ref{eq:dist1} and \ref{eq:dist2} in Section \ref{sec:pointcloud}, we see that the distances $|y_i|$ and $\sqrt{y_i^2+y_j^2}$ are not sensitive to signs (positive/negative) of the components $y_i$. For instance, there would be difficulty in distinguishing between say, $X'=(1,2)$ and $Y'=(-1,2)$. Thus, we can see that our proposed method works better if components in the data objects are all positive (or mostly positive). 

After the standardization step, each component $x_i$ of the data object $X=(x_1,x_2,\dots,x_m)$ comes from a distribution with mean 0 and standard deviation 1. Hence, by adding the fixed vector $\mathbf{v}=(5,6,7,\dots,m+4)$ (note that all components of $\mathbf{v}$ are 5 and above) to $X$, we have taken reasonable steps to try to make components of $X'=X+\mathbf{v}$ mostly positive, since only components that are more than 5 standard deviations below the mean (in $X$) would remain negative in $X'$. We remark that the above choice of fixed vector $\mathbf{v}$ is not unique (there could be other choices of $\mathbf{v}$ that work as well).

\section{Experimental Results}
\label{sec:experiment}
To evaluate the effectiveness of our proposed method, we use a real-world mixed dataset on heart disease\cite{detrano1989international} taken from the UCI Machine Learning Repository.\cite{Dua:2019} We focus on the Cleveland dataset, which comprises of data from patients referred for coronary angiography at the Cleveland Clinic. The algorithms were mostly implemented in Python, with the exception of computing persistence diagrams and Wasserstein distances using the R package \texttt{TDA}.\cite{fasy2014introduction} The codes in the paper are made publicly available on GitHub: \texttt{https://github.com/wuchengyuan88/topology-mixed-data}.

The 14 attributes of the heart disease dataset along with their data types and a brief description are presented in Table \ref{table:attributes}. For this dataset, heart disease is defined as greater than $50\%$ luminal narrowing of any major epicardial vessel.

\begin{table}[!htp]
\tbl{Attributes of the heart disease dataset (6 numeric and 7 categorical predictive attributes, and 1 categorical target attribute).}
{
\footnotesize
\begin{tabular}{l l l l}
\hline
   & Attribute & Data type        & Brief description\\
\hline
1  & age       & numeric     & Age in years                                                                             \\
2  & sex       & categorical & Gender of patient (1 = male; 0 = female)                                         \\
3  & cp        & categorical & Chest pain type (1 = typical angina; 2 = atypical angina; \\
&&&3 = non-anginal pain; 4 = asymptomatic)\\
4  & trestbps  & numeric            &  Resting blood pressure in mmHg  \\
5  & chol      & numeric            & Serum cholesterol in mg/dl     \\
6  & fbs       &  categorical           & Fasting blood sugar $>$ 120 mg/dl (1 = true; 0 = false)  \\
7  & restecg   & categorical             & Resting electrocardiographic results (0 = normal; \\
&&& 1 = having ST-T wave abnormality; \\
&&& 2 = left ventricular hypertrophy)  \\
8  & thalach   & numeric            & Maximum heart rate achieved \\
9  & exang     & categorical            & Exercise induced angina (1 = yes; 0 = no) \\
10 & oldpeak   &numeric             & ST depression induced by exercise relative to rest \\
11 & slope     &categorical             & Slope of the peak exercise ST segment (1 = upsloping;\\
&&& 2 = flat; 3 = downsloping) \\
12 & ca        & numeric            & Number of major vessels (0-3) colored by fluoroscopy \\
13 & thal      & categorical            & Heart status (3 = normal; 6 = fixed defect; \\
&&& 7 = reversible defect) \\
14 & num (target attribute)      & categorical            &  Presence of heart disease (0 = healthy; 1 = heart disease)\\
\hline                     
\end{tabular}}
\label{table:attributes}
\end{table}

After one-hot encoding, there are a total of 25 predictive attributes. Hence, each patient is represented as a data object $X$ in $\R^m$, where $m=25$. Each patient will then be represented by a point cloud $S(X')$ consisting of $25+1=26$ points in $\R^{25}$.

\begin{figure}[!htbp]
\begin{center}
\includegraphics[scale=0.5]{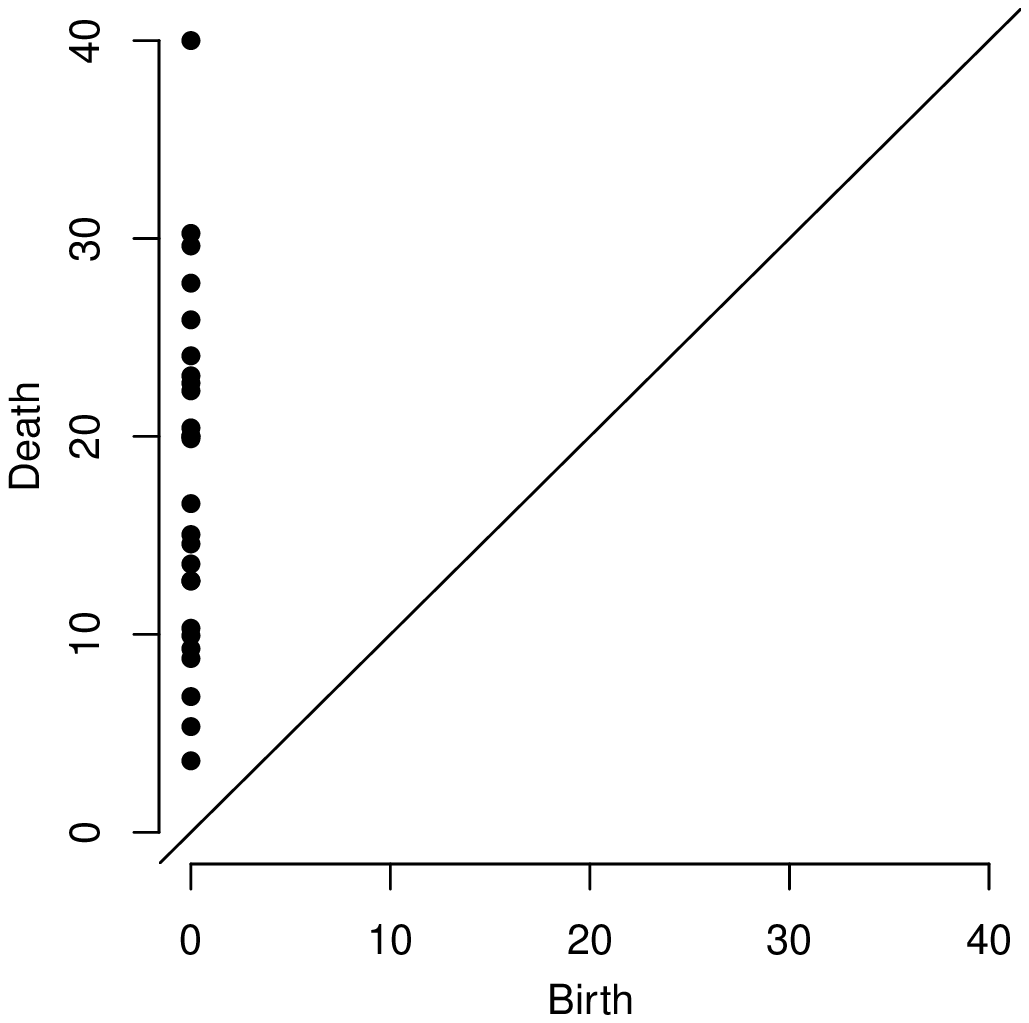}
\includegraphics[scale=0.5]{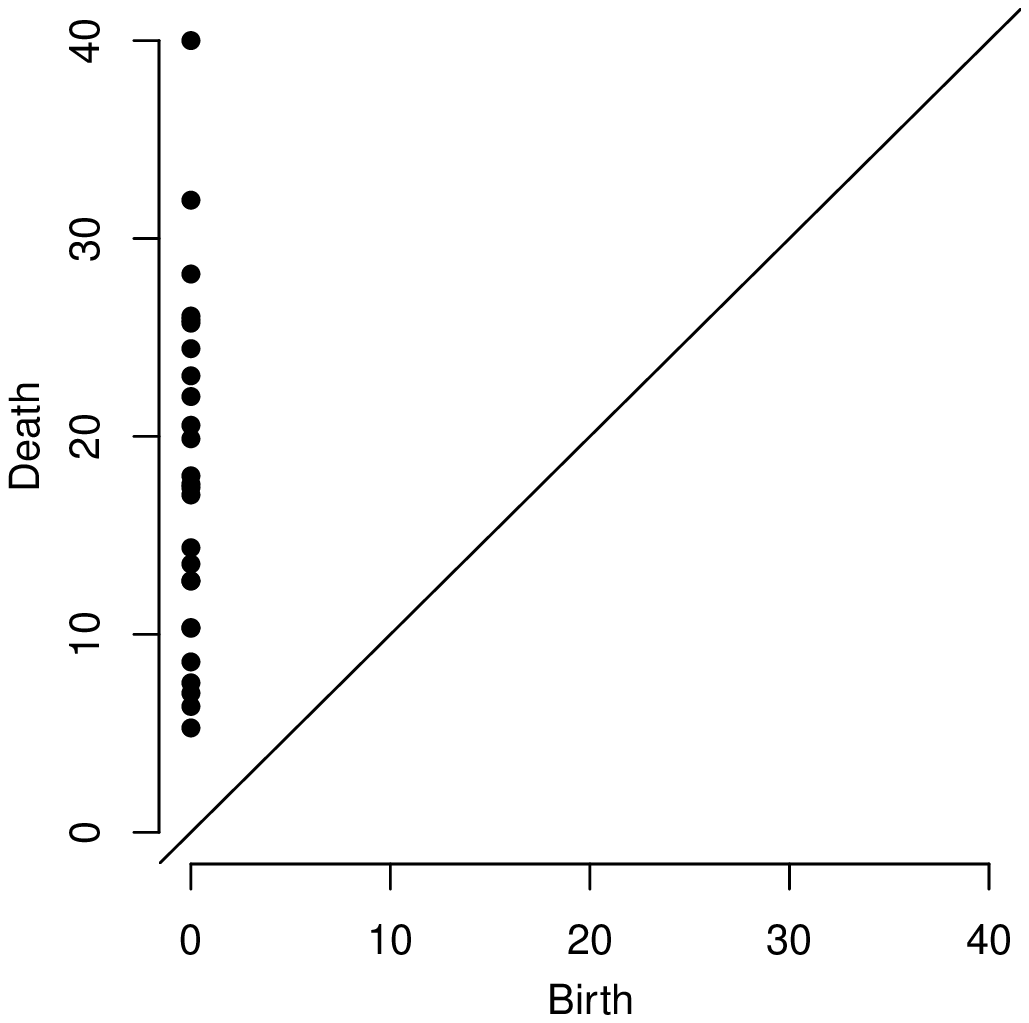}
\caption{The persistence diagram on the left belongs to a patient of class 0 (healthy), while that on the right belongs to a patient of class 1 (heart disease). The points refer to homological features in dimension 0.}
\label{fig:twopd}
\end{center}
\end{figure}

For the construction of persistence diagrams, we use the \texttt{ripsDiag} function in the R package \texttt{TDA}. We show examples of two persistence diagrams from different classes in Figure \ref{fig:twopd}. Qualitatively, we can visually observe some differences, for instance the persistence diagram for the patient of class 1 (heart disease) contains a denser cluster of points in the region corresponding to low death times of approximately 5 to 8. Quantitatively, the difference between persistence diagrams is measured by the Wasserstein distance, using the \texttt{wasserstein} function from the R package \texttt{TDA}. For this paper, distances between persistence diagrams are computed using 0 dimensional features, as we experimentally observe that 1 dimensional and higher features rarely appear in the persistence diagrams for our dataset.

We split our initial dataset (consisting of 297 patients in the Cleveland heart disease dataset) randomly into training, validation and test sets in a 60:20:20 ratio. A further summary of the split data sets can be found in Table \ref{table:datasetsummary}.

\begin{table}[!htp]
\tbl{Description of split data sets.}
{
\footnotesize
\begin{tabular}{llll}
\hline
             &                        & \multicolumn{2}{l}{Data class distribution (\%)} \\\cline{3-4}
Data set     & Number of patients & 0 (healthy)          & 1 (heart disease)         \\ 
\hline
Training set & 179                    & 54.19                      & 45.81                 \\
Validation set   & 59                    & 52.54                      & 47.46                 \\
Test set   & 59                    & 54.24                        & 45.76                   \\
\hline
\end{tabular}}
\label{table:datasetsummary}
\end{table}

To choose a suitable value for the parameter $k$ in the $k$-NN algorithm, we experiment with various values of $k$ on the validation set. The $k$ nearest neighbors will be selected from the training set based on the Wasserstein distance. We show the accuracy, sensitivity (true positive rate) and specificity (true negative rate) for various values of $k$ in Table \ref{table:differentk}. We select $k=5$ as it corresponds to the highest accuracy, as well as relatively high sensitivity and specificity (above $70\%$).

\begin{table}[!htp]
\tbl{Accuracy, sensitivity and specificity for different values of $k$ on the validation set.}
{\footnotesize
\begin{tabular}{lllllllllll}
\hline
Value of $k$     & 1 & 2 & 3 & 4 & \textbf{5} & 6 & 7 & 8 & 9 & 10 \\
\hline
Accuracy (\%)    & 69.49 & 74.58 & 76.27 & 77.97 & \textbf{81.36} & 79.66 & 77.97 & 76.27 & 77.97 & 79.66  \\
Sensitivity (\%) & 64.29 & 85.71 & 78.57 & 82.14 & \textbf{75.00} & 78.57 & 71.43 & 71.43 & 67.86 & 75.00  \\
Specificity (\%) & 74.19 & 64.52 & 74.19 & 74.19 & \textbf{87.10} & 80.65 & 83.87 & 80.65 & 87.10 & 83.87 \\
\hline
\end{tabular}}
\label{table:differentk}
\end{table}

With the chosen value of $k=5$, we show the results for the test set in Table \ref{table:finalresult}. We achieve a high level of accuracy, sensitivity (recall of positive class) and specificity (recall of negative class) on the test set. 

\begin{table}[!htp]
\tbl{Results for test set (using $k=5$).}
{
\footnotesize
\begin{tabular}{lllllll}
\hline
&               &             & \multicolumn{2}{c}{Precision (\%)} & \multicolumn{2}{c}{$F_1$ score (\%)} \\\cline{4-5}\cline{6-7}
Accuracy (\%) & Sensitivity (\%) & Specificity (\%) & (class 0)     & (class 1)     & (class 0)      & (class 1)      \\
\hline
\textbf{89.83}     & 88.89        & 90.62       & 90.62         & 88.89          & 90.62          & 88.89           \\
\hline
\end{tabular}}
\label{table:finalresult}
\end{table}

Following best practices in machine learning, we also report the results for 10-fold cross-validation. The optimal value of $k$ (for the $k$-NN algorithm) in the case of 10-fold cross-validation is found to be $k=16$. We list the results in Table \ref{table:10fold}.

\begin{table}[!htp]
\tbl{Results for 10-fold cross-validation (using $k=16$).}
{
\footnotesize
\begin{tabular}{lllllll}
\hline
&               &             & \multicolumn{2}{c}{Precision (\%)} & \multicolumn{2}{c}{$F_1$ score (\%)} \\\cline{4-5}\cline{6-7}
Accuracy (\%) & Sensitivity (\%) & Specificity (\%) & (class 0)     & (class 1)     & (class 0)      & (class 1)      \\
\hline
82.52     & 79.51        & 85.54      & 82.89         & 82.10         & 83.90          & 80.37         \\
\hline
\end{tabular}}
\label{table:10fold}
\end{table}

For reference, the accuracy of state-of-the-art algorithms reported in the literature typically ranges from around $60\%$ to $90\%$.\cite{cheung2001machine,das2009effective,latha2019improving,nahar2013computational,polat2007automatic,pouriyeh2017comprehensive} We compare our test set results with some previous results reported in literature. We remark that some of  the accuracy results in the literature are based on 10-fold cross-validation on the total data, while other results are based on train-test splits. The 10-fold cross-validation methodology has the advantage of reduced bias as every data point gets to be tested exactly once and is used in training 9 times. However, some authors also argued that selecting the best training parameters on a validation set and reporting prediction on a test set (which is how we obtained our test accuracy) is more authentic than simply performing a 10-fold cross-validation on a training set.\cite{nahar2013computational} 

Table \ref{table:comparison} gives the classification accuracies of our method and other previous approaches. Our Topological Machine Learning for Mixed Data method (TopMix) outperforms a number of other algorithms, including several state-of-the-art algorithms.

We remark that the Cleveland dataset in the UCI Machine Learning Repository consists of 303 original instances (including 297 complete instances and 6 instances with missing attributes). We only use the 297 complete instances (approximately 98\% of the full dataset) for confirming the efficiency of our method. For the results listed in Table \ref{table:comparison}, it is not clearly specified whether the authors used the reduced dataset of 297 instances or the full dataset. Hence, we also reprogram some of the methods in the list and show their accuracy for a fairer comparison. We use the Scikit-learn package in Python and reprogram 5 methods (SVM, Logistic regression, Decision tree, Naive Bayes, and Multi-layer Perceptron) using the reduced dataset. The results are recorded in Table \ref{table:comparison} as well.

\begin{table}[!htp]
\tbl{Classification accuracies obtained with our proposed Topological Machine Learning for Mixed Data method ({TopMix}) and other classifiers from literature. We have also reprogrammed some of the methods ourselves using classifiers from Scikit-learn, labelled as Scikit-learn (2020) under the ``Source'' column.}
{
\scriptsize
\begin{tabular}{lll}
\hline
Source & Method & Accuracy (\%)\\
\hline
ToolDiag &IB1-4 &50.00\\
WEKA, RA & InductH &58.50\\
ToolDiag, RA & RBF &60.00\\
WEKA, RA &FOIL &64.00\\
ToolDiag, RA &MLP+BP &65.60\\
Scikit-learn (2020) &Decision tree &67.80\\
WEKA, RA &T2 &68.10\\
S.\ Pouriyeh et al.\ (2017) &SCRL &69.96\\
WEKA, RA &1R &71.40\\
WEKA, RA &IB1c &74.00\\
Scikit-learn (2020) &Naive Bayes &76.27\\
J.\ Nahar et al.\ (2013) &J48 &76.57\\
WEKA, RA &K\textsuperscript{*} &76.70\\
J.\ Nahar et al.\ (2013) &IBK &76.90\\
R.\ Detrano &Logistic regression &77.00\\
S.\ Pouriyeh et al.\ (2017) &Decision tree &77.55\\
J.\ Nahar et al.\ (2013) &AdaBoostM1+CFS &77.94\\
J.\ Nahar et al.\ (2013) &SMO+MFS &77.95\\
N. Cheung (2001) &BNNF &80.96\\
N. Cheung (2001) &BNND &81.11\\
N. Cheung (2001) &C4.5 &81.11\\
N. Cheung (2001) &Naive Bayes &81.48\\
J.\ Nahar et al.\ (2013) &PART &81.52\\
S.\ Pouriyeh et al.\ (2017) &SVM &84.15\\
Scikit-learn (2020) &Multi-layer Perceptron &84.75\\
J.\ Nahar et al.\ (2013) &PART+MFS &86.77\\
Polat et al.\ (2006) &Fuzzy-AIRS-Knn based system &87.00\\
Scikit-learn (2020) &Logistic regression &88.14\\
Scikit-learn (2020) &SVM &88.14\\
R.\ Das et al.\ (2009) & Neural networks ensemble & 89.01\\
\textbf{TopMix (Test accuracy)} & \textbf{Topological machine learning} &\textbf{89.83}\\
\hline
\end{tabular}}
\label{table:comparison}
\end{table}

\section{Conclusions}
\label{sec:conclusion}
Data objects with mixed numeric and categorical attributes are common in real-world applications. However, many algorithms are not compatible with mixed data and can only work on single-type data, that is, either numeric or categorical data. On the other hand, TDA is a rapidly emerging machine learning method that has benefits of robustness to noise and effectiveness in high dimensions. However, traditionally TDA is applied to point cloud data, not mixed data.

In this paper, we proposed a novel topological machine learning method to classify mixed numeric and categorical data. In our method, we utilize theory from TDA such as persistent homology, persistence diagrams and Wasserstein distance in order to study mixed data. In doing so, we expand the repertoire of TDA to include mixed data.

We test our proposed method on a heart disease dataset from the UCI machine learning repository. The experimental results demonstrate that the proposed method is effective at predicting heart disease, and also outperforms several state-of-the-art algorithms.

In conclusion, our paper represents a first step towards using TDA to classify mixed numeric and categorical data and can be viewed as a proof of concept that methods from TDA are effective in the domain of mixed data, as well as in heart disease prediction.

\section*{Acknowledgments}
The authors wish to thank the referees most warmly for numerous suggestions that have improved the exposition of this paper.



\bibliographystyle{ws-ijait}
\bibliography{bibtex12}
\end{document}